\documentclass[12pt]{amsart}

\usepackage{ljm-auth}

\author{N.I.\,Zhukova$^{\dag}$, A.Yu.\,Dolgonosova$^{\ddag}$}
\crauthor{N.I.\,Zhukova, A.Yu.\,Dolgonosova} 

\tit{Pseudo-Riemannian foliations and their graphs}
\shorttit{Pseudo-Riemannian foliations and their graphs} 

\usepackage{graphicx}
\usepackage{amsfonts}
\usepackage{mathptmx}
\usepackage{amsmath}

\usepackage{amsopn, amsfonts, euscript,amscd,amssymb,latexsym}

\newtheorem{definition}{Definition}
\newtheorem{theorem}{Theorem}
\newtheorem{corollary}{Corollary}
\newtheorem{lemma}{Lemma}
\newtheorem{remark}{Remark}
\newtheorem{proposition}{Proposition}
\begin{document}
\address{}

\date{Received: date / Accepted: date}

\maketit

\address{{$\dag$}, {$\ddag$} National Research University Higher School of Economics \\
Bolshaja Pecherskaja str., 25/12, Nizhny Novgorod,  603155, Russia
}

\email{nina.i.zhukova@gmail.com, annadolgonosova@gmail.com}

\abstract{We prove that a foliation $(M, F)$ of codimension $q$ on a $n$-dimen\-sio\-nal pseudo-Riemannian manifold
	is pseudo-Riemannian if and only if any geodesic that is orthogonal at one point
	to a leaf is orthogonal to every leaf it intersects.
	
	We show that on the graph $G = G(F)$ of a pseudo-Riemannian foliation
	there exists a unique pseudo-Riemannian metric such that canonical projections are
	pse\-u\-do-Rieman\-ni\-an submersions and the fibres of different projections are 
	orthogonal at common points.
	Relatively this metric the induced foliation $(G,\mathbb{F})$ on the graph is 
	pseudo-Riemannian and the structure of the leaves of $(G,\mathbb{F})$ is described. 
	Special attention is given 	to the structure of graphs of transversally (geodesically) 
	complete pseudo-Riemannian foliations and totally geodesic pseudo-Riemannian ones.} 
	\notes{0}{
\thank{Partially supported by the Russian Foundation of Basic Research (grant no. 16-01-00132)
and by the Program of Basic Research Program at the National Research University Higher School of 
Economics in 2016 (project no. 98)}
}

\section{Introduction}
\label{intro}

Let $(M, F)$ be a smooth foliation. Recall that pseudo-Riemannian metric $g$ on the manifold $M$ 
is transversally projectable if the Lie derivative $L_{X}g$ along $X$ is zero for any vector 
field $X$ tangent to this foliation.

\begin{definition}\label{D1}
	A foliation on a pseudo-Riemannian manifold $(M, F)$ is referred to as a pseudo-Riemannian foliation if every leaf $L$ with induced metric is a pseudo-Riemannian manifold and $g$ is transversally projectable.
\end{definition}

A pseudo-Riemannian submersion (see \cite{O'Neill}) is a smooth map $p: M\to B$ which is onto and 
satisfies the following three axioms: 

(a) the differential $p_{*x} : T_xM \to T_{p(x)}B$ is onto for all $x\in M$;

(b) the fibres $p^{-1}(b)$, $b\in B$, are pseudo-Riemannian submanifolds of $M$;

(c) the differential $p_*$ preserves scalar products of vectors normal to fibres.

Definition \ref{D1} is equivalent to the fact that foliation $(M, F)$ is given locally 
by pseudo-Riemannian submersions, i.e. it is equivalent to the following definition. 

\begin{definition}\label{D2}
	A foliation $(M, F)$ on a pseudo-Riemannian manifold $(M, g)$ is said to be a pseudo-Riemannian 
	foliation if at any point there exists an adapted neighborhood $U$ and a Riemannian metric 
	$g^V$ on the leaf space $V=U/F_{U}$ such that the canonical projection $f: U\to V$ 
	is a pseudo-Riemannian submersion $(U, g|_{U})$ onto $(V, g^V)$.
\end{definition}

Let an operator of covariant differentiation be given on a manifold $M$. It is equivalent to the existence 
an linear connection on $M$ and also to the existence of an affine connection on $M$ \cite{KN}.
It is said to be {\it a manifold of linear connection.} 

Pseudo-Riemannian foliations include Lorenzian and Riemannian foliations and may be considered as 
foliations with transverse linear connection. In the foliation category isomorphisms of two foliations 
are diffeomorphisms transforming leaves of one to the leaves of the other. We proved 
[\cite{ZhD},Theo\-rem~1.1] that the group of all automorphisms of a foliation with transverse 
linear connection admits a structure of an infinite dimensional Lie group modelled on $LF$-spaces. 
Therefore the same statement is valid also for every pseudo-Riemannian foliation.

Further the pseudo-Rieman\-ni\-an manifolds $(M,g)$ are considered
with the Levi-Civita connection
$\nabla$.

\begin{definition} Let $(M,\nabla)$ be a manifold $M$ with a linear connection 
	$\nabla$. A smooth distribution $D$ on a manifold $M$
	is called {\it geodesically invariant} if for any point $x\in M$
	and each vector $X\in D_x$ the geodesic $\gamma = \gamma(s)$ of $(M,\nabla)$
	such that $\gamma(0) = x$ and $\dot{\gamma}(0) = X$ has the property
	$\dot{\gamma}(s) \in D_{\gamma(s)}$ for each $s$ of the domain of $\gamma$.
	
	A foliation $(M,F)$ on $(M,\nabla)$ is called geodesically invariant or 
	totally geodesic if its tangent distribution $TF$ is geodesically invariant. 
\end{definition}

First we prove the following criterion of a pseudo-Riemannian nature for a foliation
of a pseudo-Riemannian manifold.

\begin{theorem}\label{t1} Let $(M,F)$ be a foliation of codimension $q$ on an 
	$n$-dimensional pseudo-Rieman\-ni\-an manifold $(M,g),$ $0 < q < n$.
	Then $(M,F)$ is a pseudo-Riemannian foliation if and only if the
	$q$-dimensional distribution $D$, orthogonal to $TF$, is geodesically invariant 
	and the metric on the leaves is non-degenerate.
\end{theorem}

\begin{corollary}
	Let $M$ and $B$ be pseudo-Riemannian manifolds and $p: M \to B$ be surjective submersion. 
	Then this submersion is pseudo-Riemannian iff there is induced a pseudo-Riemannian metric 
	on the fibers and any geodesic orthogonal to the fiber at one point also is orthogonal 
	to every fiber it intersects.
\end{corollary}

For Riemannian foliations a similar result was proposed by B. Reinhardt
\cite{Rei}, and it was proven in detail by P. Molino [\cite{Mo},
Propositions 3.5 and 6.1].
We emphasize that Molino's proof bear on the property of a geodesic to be
a local extremum of the length  functional that does not have
an analogue in pseudo-Riemannian geometry.

In the proof of Theorem \ref{t1} we essentially use
the result of A.D. Lewis~\cite{Lev} on geodesically invariance
of distributions on manifolds with affine connection.

By a holonomy group $\Gamma(L)$ of a leaf $L$ of a foliation $(M, F)$ we mean a germ holonomy group of $L$ usually used in the foliation theory \cite{Tam}. If $\Gamma(L)=0$ the leaf $L$ is said to be a leaf without holonomy.

A construction of the holonomy groupoid of a foliation was presented by S.~Ehresmann.
Another equivalent construction was given
by H. Winkelnkemper~\cite{Win}
and named by him {\it the graph} of a foliation. The graph  $G(F)$ contains all
information about the foliation $(M, F)$ and its holonomy groups. $C^*$-algebras
of complex valued functions of foliations $(M, F)$ are determined on
$G(F)$  and are one of the fundamental concepts in $K$-theory of
foliations.

In the general case the graph of a smooth foliation $(M, F)$ of codimension
$q$ on an $n$-di\-men\-si\-onal smooth manifold $M$ is a non-Hausdorff smooth
$(2n-q)$-manifold (the precise definition is given in Subsection \ref{sec:2.1}).

For the graph $G(F)$ of a  pseudo-Riemannian  foliation $(M,F)$ we prove the following statement.

\begin{theorem}\label{t2}
	Let $(M,F)$ be a smooth pseudo-Riemannian foliation of
	codimension  $q$ on the $n$-dimensional pseudo-Riemannian manifold $(M,g)$.
	Let $G(F)$ be its graph with the canonical projections
	$p_i : G(F) \to M$, $i=1,2$.
	Then:
	
	\begin{enumerate}
		\item The graph $G(F)$ of a foliation $(M,F)$  is a Hausdorff $(2n-q)$-dimen\-sio\-nal 
		manifold with the induced foliation 
		$\mathbb F = \left\lbrace \mathbb L_\alpha = p^{-1}_i(L_{\alpha}) | L_{\alpha}\in{F}\right\rbrace$, 
		$i = 1,2.$ Moreover the germ holonomy groups $\Gamma({L}_\alpha)$ and
		$\Gamma(\mathbb{L}_\alpha)$ of the appropriate leaves
		$L_\alpha$ and $\mathbb{L}_\alpha$ are isomorphic.
		
		\item On the graph $G({F})$ there exists a unique pseudo-Riemannian metric $d$
		with respect to which $(G({F}),\mathbb F)$ is a pseudo-Riemannian foliation
		and $p_i$ are pseudo-Riemannian submersions. In this case fibres of $p_1$ are orthogonal 
		to the fibres of $p_2$ at common points.
		
		\item Every leaf $\mathbb{L}_\alpha = p^{-1}_i(L_\alpha)\in\mathbb F$ is a reducible pseudo-Riemannian
		ma\-nifold that is isometric to the quotient manifold
		$({\mathcal L}_\alpha\times {\mathcal L}_\alpha)/{\Psi}_{\alpha}$ of the
		pseudo-Riemannian product ${\mathcal L}_\alpha\times {\mathcal L}_\alpha$
		of the pseudo-Riemannian holonomy covering space ${\mathcal L}_\alpha$ of $L_\alpha$
		by the isometry group $\Psi_{\alpha}$, and $\Psi_{\alpha}\cong\Gamma(L_\alpha)\cong\Gamma(\mathbb{L}_\alpha)$.
	\end{enumerate}
\end{theorem}
\begin{definition} The pseudo-Riemannian metric $d$ on the graph $G(F)$
	satisfying Theorem \ref{t2} is called the {\it induced metric}.
\end{definition}
\begin{corollary} There exists a dense saturated $G_\delta$-subset of $G(F)$ any leaf $(\mathbb{L}_\alpha, d)$
	of which is isometric to the direct product $L_\alpha\times L_\alpha$ of a pseudo-Riemannian manifolds $(L_{\alpha}, g).$
\end{corollary}

\begin{remark}
	We prove property $3$ in Theorem \ref{t2} without the application of
	the well known Wu's theorem \cite{Wu}.
	This theorem is not applicable here because
	the completeness of the pseudo-Riemannian metric $d$ is not assumed.
\end{remark}

By a geodesic we mean a piecewise geodesic.
\begin{definition}
	A pseudo-Riemannian foliation $(M,F)$ is called {\it transversally complete} if the canonical parameter on every maximal orthogonal geodesic is defined on the whole real line.
\end{definition}
Under the additional assumption of the transversal completeness
of a pseudo-Rie\-man\-ni\-an foliation we prove the following statement.

\begin{theorem}\label{t3}
	Let $(M,F)$ be a transversally complete pseudo-Riemannian foliation on a pseudo-Riemannian manifold
	$(M,g)$ and  $d$ be the induced pseudo-Riemannian
	metric on its graph $G(F)$. Then I:
	\begin{enumerate}
		\item The induced foliation $(G(F),\mathbb F)$ is a transversally complete pseudo-Rie\-man\-ni\-an foliation.
		
		\item The orthogonal $q$-dimensional distributions $\mathfrak M$
		and $\mathfrak N$ are Ehresmann connections for the foliations $(M,F)$
		and $(G(F),\mathbb F)$ respectively, and, for any $L_{\alpha}\in{F}$,
		$\mathbb{L_\alpha} = p^{-1}_i(L_\alpha)$, the following holonomy groups $\Gamma(L_\alpha)$,
		$H_{\mathfrak M}(L_\alpha)$, $\Gamma(\mathbb{L_\alpha})$ and 
		$H_{\mathfrak N}(\mathbb{L}_{\alpha})$ are isomorphic.
		
		\item The canonical projections $p_i: G(F) \to M$, $i=1,2$, form
		locally trivial fibrations with the same standard fibre $L_0$,
		and $L_0$ is diffeomorphic to any leaf without holonomy of $(M,F)$.
		
		\item Every leaf $\mathbb{L}_\alpha = p^{-1}_i(L_\alpha)\in\mathbb F$
		is a reducible pseudo-Riemannian manifold that is isometric
		to the quotient manifold
		$({\mathcal L}_\alpha\times {\mathcal L}_\alpha) /{\Psi_{\alpha}}$ of the
		pseudo-Riemannian product ${\mathcal L}_\alpha\times {\mathcal L}_\alpha$
		of the pseudo-Riemannian holonomy covering space ${\mathcal L}_\alpha$
		for $L_\alpha$ by an isometry group $\Psi_{\alpha}$, where
		$\Psi_{\alpha}\cong\Gamma(L_\alpha)\cong H_{\mathfrak M}(L_\alpha)\cong
		\Gamma(\mathbb{L}_{\alpha})\cong H_{\mathfrak N}(\mathbb{L}_{\alpha})$,
		and every ${\mathcal L}_\alpha$ is diffeomorphic to $L_0$.
	\end{enumerate}
	
	II. If moreover, the foliation $(M,F)$ is also geodesically invariant, then:
	\begin{enumerate}
		\item[(i)] Each foliation $\mathbb{F},F^{(i)}:=\left\lbrace p^{-1}_{i}(x)|x\in M \right\rbrace,\, i=1,2, $ is geodesically invariant and pseudo-Riemannian.
		\item[(ii)] Any leaf without of holonomy of $(M, F)$ is isometric to any fibres of submersions 
		$p_1$ and $p_2$ with respect to corresponding induced metrics.
		\item[(iii)] Every leaf $\mathbb {L}$ without holonomy of the foliation $(G(F),\mathbb{F})$ 
		is isometric to the pseudo-Riemannian product $L_0\times L_0$ and any other leaf $\mathbb{L}_\alpha$ 
		is isometric to the pseudo-Riemannian quotient manifold $(L_0\times L_0)/\Psi_{\alpha}$, where
		$\Psi_{\alpha}\cong\Gamma(\mathbb{L}_{\alpha})$.
	\end{enumerate}
\end{theorem}

We emphasize that the proofs of Theorems \ref{t2} and \ref{t3} essentially use
the notion of an Ehresmann connection of a foliation proposed by R.A.~Blumenthal and
J.J~Hebda~\cite{BH} and also the results of the first author on graphs of foliations with Ehresmann connection \cite{ZhG}, \cite{ZhL} and some other statements (see Sections \ref{sec:4}-\ref{sec:5}).

\medskip{\noindent\bf Notations\,} We denote by ${\mathfrak X}(N)$ the set of all
smooth vector fields on a manifold $N.$ If $\mathfrak M$ is a smooth
distribution on $M$, then ${\mathfrak X}_{\mathfrak M}(M):=\{X\in{\mathfrak X}(M)\mid
X_u\in {\mathfrak M}_u\,\,\,\,\,\,\forall u\in M\}$. For a foliation $(M,F)$
we denote ${\mathfrak X}_{TF}(M)$ also by ${\mathfrak X}_{F}(M)$.
Let us denote the leaf of foliation $(M, F)$ passing through a point $x\in M$ by $L(x).$

Let $\mathfrak F\mathfrak o\mathfrak l$ be the category of foliations where morphisms are
smooth maps transforming leaves into leaves.

Let $\cong$ be the symbol of a group isomorphism and also of
a manifold diffeomorphism.

Further smoothness is understood to mean $C^{\infty}.$

\section{A criterion of pseudo-Riemanniance of a foliation}
\label{sec:2}
\subsection{Foliate and transversal vector fields}
Let $(M, F)$ be a smooth foliation of codimension $q$ of a smooth $n$-dimensional manifold $M.$

A function $f\in{\mathfrak{F}}(M)$ is called {\it basic} if it is constant on 
every leaf $L$ of this foliation. 	Let $\Omega^{0}_b(M,F)$ be a set of basic functions. 
Remark that $f\in\Omega^{0}_b(M,F)$ if and only if $f$ depends on only 
transverse coordinates in any adapted charts.
\begin{definition}\label{d3}
	Let $(M,F)$ be an arbitrary foliation. A vector field $X\in\mathfrak{X}(M)$
	is called {\it foliate} if for any $Y\in\mathfrak{X}_{F}(M)$ the vector field $[X, Y]$ 
	belongs to $\mathfrak{X}_{F}(M)$.
\end{definition}
Following Molino \cite{Mo} we denote by $L(M, F)$ the set of all foliate vector fields.
In this case $L(M,F)$ is a normalizer in $\mathfrak{X}(M)$ of the Lie subalgebra $\mathfrak{X}_{F}(M).$
Therefore $L(M, F)$ is a Lie subalgebra of $\mathfrak{X}(M)$.

A $q$-dimensional smooth distribution $\mathfrak M$ on the manifold $M$ is called 
{\it trans\-versal} to the foliation $(M, F)$ if the equality $T_{x}M= T_x{F}\oplus{\mathfrak M}_x$ 
holds for any $x\in M$, where ${\oplus}$ stands for the direct sum of 
vector spaces. Let us identify the quotient vector bundle $TM/TF$ with 
a transversal distribution $\mathfrak{M}$.

Let $K_{a}^{m}$ be an open cube in $R^{m}$ with the center at the null $0_{m}\in R^{m}$ 
and with the side $2a.$ 	A chart $(U,\varphi)$ of the manifold $M$ is said to be {\it adapted}
with respect to the foliation $(M, F)$ at $p_{0}\in U$ if $\varphi(U) = K_{a}^{n-q}\times K_{a}^{q}$,
and $\varphi(p_{0})=(0_{n-q},0_{q}),\,\varphi(p)=(x^1 ,...,x^{n-q},y^{1},...,y^{q}) =
(x^{a},y^{\alpha})\in K_{a}^{n-q}\times K_{a}^{q}$ and 
$\varphi^{-1}(K_{a}^{n-q}\times\{y^{\alpha}\}), \{y^{\alpha}\}\in K_{a}^{q}$, 
is a plaque of $p$ in $U,$ i.e. it is a connected component of $p$ in $L\cap U.$ 
We say that $x^{a}$ are leaf coordinates and $y^{\alpha}$ are transverse coordinates.
At every point $x\in U$ the derivations  
$\frac{\partial}{\partial y^1},...,\frac{\partial}{\partial y^q}$ 
generate a subspace of $T_xM$ that is transversal to $T_xF$. Vectors 
$\frac{\overline{\partial}}{\partial y^1},...,\frac{\overline{\partial}}{\partial y^q}$
or for short $\frac{\overline{\partial}}{\partial y^\beta},\,\, \beta=1,...,q$, are
the corresponding local sections of the distribution $\mathfrak{M}_{x}$ at each point
$x\in U.$ The section $\overline{X}$ of $\mathfrak{M}$ determined by
the vector field $X\in\mathfrak{X}(M)$ has in $U$ the form
$\overline{X} = \sum_{\beta=1}^{q}\overline{X}^{\beta}\frac{\overline{\partial}}{\partial y^\beta},$
where $\overline{X}^{i}$ are the last $q$ components of ${X}$ in the local basis 
$\{\frac{\partial}{\partial x^a},\frac{\partial}{\partial y^\beta}\}$. 

In particular, if $X$ is a foliate vector field, then $\overline{X}^{\beta}$ depend only on
the transverse 	variables $y^1,...,y^{q}$. In this case the vector field $\overline{X}$ 
is called the {\it transverse vector field} associated to $X.$ Let $l(M, F)$ be the set of
transverse vector fields. Then the projection $L(M,F)\to l(M, F): X\mapsto\overline{X}$
is well defined with kernel is equal to $\mathfrak{X}_{F}(M).$ Therefore, for the foliation 
$(M, F)$ there exists an exact sequence of vector spaces
\begin{equation}
0 \to{\mathfrak X}_{F}(M) \to L(M,F) \to l(M,F) \to 0
\end{equation}

Further we consider coordinate charts $(U,\varphi)$ adapted to $(M, F)$ with the base vector fields
$\{\frac{\partial}{\partial x^a},\frac{\overline{\partial}}{\partial y^\beta}\},$ where
$\frac{\overline{\partial}}{\partial y^\beta}\in l(U, F_U)$. We call
$\{\frac{\partial}{\partial x^a},\frac{\overline{\partial}}{\partial y^\beta}\}$
{\it an $\mathfrak M$-adapted foliate basis} on $U$.

The following characteristic properties of foliate vector fields are well known
(see, for example, \cite{Mo}).

\begin{proposition}
	Let $(M,F)$ be a smooth foliation  of codimension $q$ and $\mathfrak M$ be
	the transverse $q$-dimensional distribution. Then the following conditions are equivalent:
	
	\begin{enumerate}
		\item The vector field $X$ is foliate in the sense of Definition \ref{d3}.
		\item For any basic function $f\in \Omega^{0}_b(M,F)$ the function $Xf$ is also basic.
		\item For any adapted chart $(U,\varphi)$ to $(M,F)$ with coordinates $(x^{a}, y^{\beta})$
		and $\mathfrak M$-adapted foliate basis $\{\frac{\partial}{\partial x^a}, \frac{\overline{\partial}}{\partial y^\beta} \}$ 
		on $U$ the vector field $X$ has a form
		$X = \sum_{a=1}^p X^{a}\frac{\partial}{\partial x^a} +
		\sum_{\delta=1} ^{q} \overline{X}^{\delta}\frac{\overline{\partial}}{\partial y^\delta},$
		and $\overline{X}^{\delta} = \overline{X}^{\delta}(y^{\beta}).$
	\end{enumerate}
\end{proposition}

\subsection{Pseudo-Riemannian foliations of the pseudo-Riemannian manifolds}
\label{ssec:2.2}
\begin{proposition}\label{p2}
	Let $(M,F)$ be a foliation of codimension $q$ on a pseudo-Riemannian manifold $(M,g)$
	and $\mathfrak M$ be the orthogonal $q$-dimensional distribution.
	The following conditions are equivalent:
	
	\begin{enumerate}
		\item The pseudo-Riemannian metric $g$ is transversally projectable.		
		\item $g(X,Y)\in \Omega^0_b(M,F)$ for any foliate vector fields $X, Y\in {\mathfrak X}(M)$ orthogonal to this foliation.
		\item The metric $g$ in any $\mathfrak M$-adapted foliate basis 
		$\{\frac{\partial}{\partial x^a}, \frac{\overline\partial}{\partial y^\beta}\}$
		corresponding to an adapted chart has the form
		$$(g_{ij}) = \left(\begin{array}{cc}
		g_{ab}(x^d, y^\beta) & 0\\
		0 & g_{\delta\varepsilon }(y^\beta)
		\end{array}\right),
		$$ where $a,b,d =\overline{1,p}; \, \,\, \beta,\delta,\varepsilon = \overline{1,q}.$
	\end{enumerate}
\end{proposition}


\subsection{Lewis's criterion}

{\it Lewis's criterion} \cite{Lev}. {\it  A smooth distribution $\mathfrak{M}$ on a manifold of affine connection $(M,\nabla)$
	is geodesically invariant if and only if $$\nabla_{X}Y + \nabla_{Y}X \in{\mathfrak X}_{\mathfrak{M}}(M)$$
	for any vector fields $X, Y$ belonging to ${\mathfrak X}_{\mathfrak{M}}(M)$.}

\subsection{Proof of Theorem \ref{t1}}
Let $(M,F)$ be a smooth foliation of codimension $q$ on
a pseudo-Riemannian manifold $(M,g)$ such that the restriction of this metric to the leaves
is non-degenerate. Suppose now that the  $q$-dimensional distribution
$\mathfrak{M}$ orthogonal to $(M,F)$ is geodesically invariant.

Let us consider any foliate vector fields $X^\mathfrak{M}, Y^\mathfrak{M} \in{\mathfrak X}_{\mathfrak{M}}(M)$
and an arbitrary vector field $Z^F\in{\mathfrak X}_{F}(M)$.
Observe that the metric $g$ is transversally projectable with respect to
the foliation $(M,F)$ if and only if
\begin{equation}\label{e1}
Z^F\cdot g(X^\mathfrak{M},Y^\mathfrak{M}) = 0.
\end{equation}

Recall that the equality $\nabla_{X} g = 0$ is equivalent to
\begin{equation}\label{e2}
Z\cdot g(X,Y) = g(\nabla_{Z}X,Y) + g(X,\nabla_{Z}Y) \,\,\, \,\,\,\forall X,Y,Z \in{\mathfrak X}(M).
\end{equation}
Putting in (\ref{e2}) $Z=Z^F\in{\mathfrak X}_F(M)$ and $X=X^\mathfrak{M}$, 
$Y=Y^\mathfrak{M}\in {\mathfrak X}_\mathfrak{M}(M)$ and
using the identity $g(X^\mathfrak{M},Z^F) = 0$ we have
\begin{equation}\label{e3}
Y^\mathfrak{M}\cdot g(X^\mathfrak{M},Z^F) = g(\nabla_{Y^\mathfrak{M}}X^\mathfrak{M},Z^F) + 
g(X^\mathfrak{M},\nabla_{Y^\mathfrak{M}}Z^F) = 0.
\end{equation}
By analogy, changing $Y^{\mathfrak{M}}$ and $X^{\mathfrak{M}}$ we get
\begin{equation}\label{e4}
X^\mathfrak{M}\cdot g(Y^\mathfrak{M},Z^F) = g(\nabla_{X^{\mathfrak{M}}}Y^{\mathfrak{M}},Z^{F}) + 
g(Y^{\mathfrak{M}},\nabla_{X^{\mathfrak{M}}} Z^{F}) = 0.
\end{equation}
Add (\ref{e3}) to (\ref{e4}), then apply the bi-linearity of the
pseudo-Riemannian metric $g$ and obtain
\begin{equation}\label{e5}
g(\nabla_{Y^\mathfrak{M}}X^\mathfrak{M} + \nabla_{X^{\mathfrak{M}}}Y^{\mathfrak{M}},Z^F) + 
g(X^\mathfrak{M},\nabla_{Y^\mathfrak{M}}Z^F) + g(Y^{\mathfrak{M}},\nabla_{X^{\mathfrak{M}}} Z^{F}) = 0.
\end{equation}
Due to the geodesically invariance of the distribution $\mathfrak{M}$ the Lewis's criterion implies
$\nabla_{X^{\mathfrak{M}}}Y^{\mathfrak{M}} + \nabla_{Y^{\mathfrak{M}}}X^{\mathfrak{M}}\in {\mathfrak{X}}_{\mathfrak{M}}(M).$
Therefore the first term in (\ref{e5}) was equal to zero. Since $\nabla$ is the Levi-Civita
connection, it is torsion free and
\begin{equation}\label{e6}
\nabla_{Y^{\mathfrak{M}}}Z^{F} = \nabla_{Z^{F}}Y^{\mathfrak{M}} + [Y^{\mathfrak{M}},Z^{F}],
\end{equation}
similarly
\begin{equation}\label{e7}
\nabla_{X^{\mathfrak{M}}}Z^{F} = \nabla_{Z^{F}}X^{\mathfrak{M}} + [X^{\mathfrak{M}},Z^{F}].
\end{equation}
Putting (\ref{e6}) and  (\ref{e7}) into (\ref{e5}) we obtain
$$
g(X^{\mathfrak{M}}, \nabla_{Z^F}Y^{\mathfrak{M}}) + g(Y^{\mathfrak{M}}, \nabla_{Z^F}X^{\mathfrak{M}}) + g(X^{\mathfrak{M}}, [Y^{\mathfrak{M}},Z^{F}]) + g(Y^{\mathfrak{M}}, [X^{\mathfrak{M}},Z^{F}]) = 0.
$$
In concordance with the choice, the vector fields $X^{\mathfrak{M}}$ and $Y^{\mathfrak{M}}$ are foliate, so
$[Z^{\mathfrak{M}},Y^{F}]$ and $[X^{\mathfrak{M}},Y^{F}]$ belong to $\mathfrak{X}_{F}(M).$
Hence the third and fourth terms in the previous equation vanish. Therefore
\begin{equation}\label{e8}
g(X^{\mathfrak{M}}, \nabla_{Z^F}Y^{\mathfrak{M}}) + g(Y^{\mathfrak{M}}, \nabla_{Z^F}X^{\mathfrak{M}}) =0.
\end{equation}

Let in (\ref{e2}) $Z=Z^F\in {\mathfrak X}_{F}(M)$,  $X=X^\mathfrak{M}$, 
$Y=Y^\mathfrak{M}\in {\mathfrak X}_{\mathfrak{M}}(M)$ and using the
relation (\ref{e8}) we have the following
\begin{equation}\label{e9}
Z^F\cdot g(X^\mathfrak{M},Y^\mathfrak{M}) = g(\nabla_{Z^F}X^\mathfrak{M},Y^\mathfrak{M}) + 
g(X^\mathfrak{M},\nabla_{Z^F}Y^\mathfrak{M}) = 0.
\end{equation}

The equality (\ref{e9}) implies the pseudo-Riemanniance of the foliation $(M,F)$.

Converse. Let $(M,F)$ be a pseudo-Riemannian foliation of codimension $q$ on
a pseudo-Riemannian manifold $(M,g),$ hence by Definition\ref{D1} the restriction of 
this metric on leaves is non-degenerate.
Denote by $\mathfrak{M}$ the orthogonal $q$-dimensional distribution to this foliation.
The pseudo-Riemanniance of the foliation $(M,F)$ implies that equalities (\ref{e9}) and
(\ref{e8}) hold for an arbitrary vector field $Z^F\in{\mathfrak X}_{F}(M)$ and any 
foliate vector fields $X^\mathfrak{M}, Y^\mathfrak{M} \in{\mathfrak X}_{\mathfrak{M}}(M).$
From (\ref{e6}) and (\ref{e7}) we find
\begin{equation}\label{e10}
\nabla_{Z^{F}}Y^{\mathfrak{M}} = \nabla_{Y^{\mathfrak{M}}}Z^{F} - [Y^{\mathfrak{M}},Z^{F}],
\end{equation}
and similarly
\begin{equation}\label{e11}
\nabla_{Z^{F}}X^{\mathfrak{M}} = \nabla_{X^{\mathfrak{M}}}Z^{F} - [X^{\mathfrak{M}},Z^{F}].
\end{equation}
Putting (\ref{e10}) and (\ref{e11}) in (\ref{e8}) and taking into account that
$[Z^{\mathfrak{M}},Y^{F}]$ and $[X^{\mathfrak{M}},Y^{F}]$ belong to $\mathfrak{X}_{F}(M)$ we get
\begin{equation}\label{e12}
g(X^\mathfrak{M},\nabla_{Y^\mathfrak{M}}Z^F) + g(Y^{\mathfrak{M}},\nabla_{X^{\mathfrak{M}}} Z^{F}) = 0.
\end{equation}
Recall that the relation (\ref{e5}) is obtained using only the condition $\nabla g = 0$
and the orthogonality of the distributions $\mathfrak{M}$ and $TF$. Therefore we may apply
(\ref{e5}). Thus we have
\begin{equation}\label{e13}
g(\nabla_{Y^\mathfrak{M}}X^\mathfrak{M} + \nabla_{X^{\mathfrak{M}}}Y^{\mathfrak{M}},Z^F) = 0
\end{equation}
for any vector field $Z^F$ tangent to the foliation $(M,F)$.
Nondegeneracy of the restriction of $g$ to the leaves of $(M,F)$
implies
\begin{equation}\label{e14}
\nabla_{Y^\mathfrak{M}}X^\mathfrak{M} + \nabla_{X^{\mathfrak{M}}}Y^{\mathfrak{M}} \in{\mathfrak X}_{\mathfrak{M}}(M).
\end{equation}
According to the Lewis's theorem mentioned above the
relation (\ref{e14}) guarantees that the distribution $\mathfrak{M}$
is geodesically invariant.

\section{Graphs of pseudo-Riemannian foliations}
\label{sec:2}
\subsection{The graph of a smooth foliation}\label{sec:2.1}
Let $(M, F)$ be a smooth foliation of codimension $q$ of the $n$-dimensional manifold $M$
and $\mathfrak{M}$ be transversal $q$-dimension distribution on $M.$
Denote by $A(x,y)$ the set of all piecewise smooth paths from $x$ to $y$ on the same leaf $L.$
Two paths $h_1,h_2 \in A(x,y)$ are called equivalent and are denoted by
$h_{1}\sim h_{2}$ if and only if  they define the same germ at the point $x$ of the
holonomy diffeomorphisms $D_{x}^{q}\to D_{y}^q$ from $D_{x}^{q}$ to $D_{y}^q$,
where $D_{x}^{q}$ and $D_{y}^q$ are  $q$-dimensional disks transversal to this foliation.

We denote by $<h>$ the equivalence class containing $h.$

\begin{definition}
	The set
	$$G(F):=\left\lbrace (x,<h>,y)\,|\, x\in M, \, y\in L(x), \, h\in A(x,y)\right\rbrace$$
	is called {\it the graph of the foliation} $(M,F)$\cite{Win}.
\end{definition}

Suppose that $h$ and $g$ are paths such that $h(1) = g(0)$.
Denote by $h\cdot g$ the product of the path $h$ and $g$, i.e. $(h\cdot g)(t) = h(2t), t\in [0, 1/2]$
and $(h\cdot g)(t) = g(2t - 1), t\in [1/2, 1]$.
The equality $$(x,<h>,v)\circ(v,<g>,y):=(x,<h\cdot g>,y)$$
where $h(1) = g(0)$, defines a partial multiplication $\circ$
in the graph $G( F).$ In relation to this operation $G(F)$ is a groupoid
named the {\it the holonomy groupoid} of the foliation $(M,F)$.

Consider any coordinate chart $(U, \varphi)$ with the center $x\in M$ adapted to the foliation 
$(M, F)$ where $\varphi(U)=K_{a}^{n-q}\times K_{a}^{q}$. Further we call the submanifold 
$D_{x}^{q}:=\varphi^{-1}(\{0_{n-q}\}\times K_{a}^{q})$ a $q$-dimensional transversal disk at $x$. 

Let $(U, \varphi)$ and $(U^{'}, \varphi^{'})$ be two adapted to $(M, F)$ charts at $x$ and 
$x^{'}\in L(x)$ accordingly. The chart  $(U^{'}, \varphi^{'})$ is referred to as {\it a subordinate } 
to $(U,\varphi)$ if $\varphi(U) = K_{a}^{n-q}\times K_{a}^{q},$ $\varphi^{'}(U^{'}) = 
K_{a^{'}}^{n-q}\times K_{a^{'}}^{q}$, where $0<a^{'}\leq a,$ and for any $c\in K_{a}^{q}$ 
the local leaves $\varphi^{-1}(K_{a}^{n-q}\times\{c\})$ and $\varphi^{'-1}(K_{a^{'}}^{n-q}\times\{c\})$ 
belong to the same leaf of the foliation $(M, F).$ For every adapted chart $(U, \varphi)$ 
at $x$ and any $y\in L(x)$ there exists a subordinated chart $(U^{'}, \varphi^{'})$ at $y$ \cite{Pal}. 
It's not difficult to show that for any points $x,y$ from a leaf $L$ there exists a pair of
adapted charts $(U, \varphi)$ and $(U^{'}, \varphi^{'})$ such that each of them is subordinated to the other.

Let $z=(x,<h>, y)$ be any point in the graph $G(F).$ Let $(U, \varphi)$ and $(U^{'}, \varphi^{'})$ 
be a pair a pair mutually subordinated charts and $D_{x}, D_{y}$ be the $q$-dimensional transversal 
disks at $x$ and $y$ accordingly. Take any point $\widetilde{x}\in D_{x},$ and denote by $\widetilde{h}$ 
the holonomy lift of $h$ into point $\widetilde{x}$ \cite{Win}. Consider any point $x^{'}$ in the local 
leaf $L_{\widetilde{x}}$ through $\widetilde{x}$ any $y^{'}$ in the local leaf $L_{\widetilde{y}}$
through $\widetilde{y}:=\widetilde{h}(1)$. Let $t_{x^{'}}$ be a path in $L_{\widetilde{x}}$ 
connected $x^{'}$ with $\widetilde{x}$ and $t_{y^{'}}$ be a path in $L_{\widetilde{y}}$ connected 
$y^{'}$ with $\widetilde{y}.$  Denote by $V_{z}$ the set of
$z^{'}=(x^{'},<t_{x^{'}}\cdot\widetilde{h}\cdot t_{y^{'}}^{-1}> y^{'})$ 
definded by the indicated above way. The set $\varSigma:=\{V_{z}\,|\,z\in G(F)\}$ 
forms a base of a topology $\Omega$ in $G(F),$ and $(G(F),\Omega)$ is non-Hausdorff in general. 

Define a map $\chi_{z}:V_{z}\to K_{a}^{n-q}\times K_{a}^{n}$ by the following equality 
\begin{center}
	$\chi_{z}(z^{'}):=(\varphi(x^{'}),pr_{2}\circ\varphi^{'}(y^{'}))\,\,\,\,\,\,\forall z^{'}=(x^{'},
	<t_{x^{'}}\cdot\widetilde{h}\cdot t_{y^{'}}^{-1}>, y^{'}),$
\end{center} 
where $pr_{2}:K_{a}^{2n-q}\cong K_{a}^{n}\times K_{a}^{n-q}\to K_{a}^{n-q}$ is the canonical projection onto $K_{a}^{n-q}.$ The set $\mathcal{A}=\{(V_{z},\chi_{z})\,|\,z\in G(F)\}$ is a smooth atlas on $G(F).$ 

Thus, the graph $G(F)$ becomes a $(2n-q)$-dimensional smooth manifold which is non-Hausdorff, in general.

\begin{definition} A pseudogroup ${\mathcal H}$ of local holonomy diffeomorphisms
	of a manifold $N$ is called {\it quasi-analytic} if  the existence of an open connected subset
	$V$ in $N$ such that $h|_{V} = id_V$ for an element $h\in\mathcal H$ implies that
	$h = id_{D(h)}$ in the whole connected domain $D(h)$ of $h$ that contains $V$.
\end{definition}

According to (\cite{ZhL}, Proposition 2), Winkelnkemper's criterion
of the property of the graph $G(F)$ to be Hausdorff
can be reformulated in the following form:

\begin{proposition}\label{pC}
	The topological space of the graph $G(F)$ of a foliation
	$(M,F)$ is Hausdorff iff the holonomy pseudogroup of this
	foliation is quasi-analytic.
\end{proposition}

The mappings
$$p_1 : G(F) \to M:(x,<h>,y)\mapsto x,\,\,\,\,\, p_2 : G(F) \to M: (x,<h>,y)\mapsto y$$
are referred to as canonical projections, and $p_1$ and $p_2$ are submersions onto $M.$

\begin{definition} A foliation
	$\mathbb F = \{\mathbb L_\alpha = p^{-1}_1(L_{\alpha}) = p^{-1}_2(L_{\alpha})\,|\,L_{\alpha}\in {F}\}$
	is defined on the graph $G(F)$ and is called the {\it induced foliation}.
\end{definition}

\subsection{An Ehresmann connection for a smooth foliation}

Recall the notion of an Ehresmann connection which was introduced by
R.A. Blumenthal and J.J. Hebda \cite{BH}. We use the term {\it a
	vertical-horizontal homotopy} introduced previously by R.~Her\-mann. All
mappings are supposed to be piecewise smooth.

Let $(M, F)$ be a foliation of an arbitrary codimension $q\geq
1$. Let $\mathfrak{M}$ be a transversal distribution on the manifold $M$, then 
for any $x\in M$ the equality $T_{x}M= T_x{F}\oplus{\mathfrak M}_x$ holds. Vectors from
$\mathfrak M_x$, $x\in M$, are called horizontal. A piecewise smooth
curve $\sigma$ is horizontal (or $\mathfrak M$-horizontal) if each of
its smooth segments is an integral curve of the distribution
$\mathfrak{M}.$ The distribution $TF$ tangent to the leaves of the
foliation $(M, F)$ is called vertical. One says that a
curve $h$ is vertical if $h$ is contained in a leaf of the
foliation $(M, F)$.

A {\it vertical-horizontal homotopy} (v.h.h. for short) is a piecewise smooth map
$H: I_1\times I_2\to M$, where $I_1=I_2=[0,1]$, such that for any $(s,t)\in I_1\times I_2$
the curve $H|_{I_1\times \{t\}}$ is horizontal and the curve $H|_{\{s\}\times I_2}$
is vertical. The pair of curves $(H|_{I_1\times\{0\}}, H|_{\{0\}\times I_2})$
is called {\it a base of the v.h.h.} $H$.
Two paths $(\sigma, h)$ with common origin $\sigma (0)=h(0)$, where $\sigma$ is a horizontal path
and $h$ is a vertical one, are called an {\it admissible pair of paths}.

A distribution $\mathfrak M$ transversal to a foliation $(M,F)$ is called an {\it
	Ehresmann connection for} $(M,F)$ if for any admissible pair of paths $(\sigma, h)$
there exists a v.h.h. with the base $(\sigma, h)$.

Let $\mathfrak M$ be an Ehresmann connection for a foliation $(M,F)$. Then  for any
admissible pair of paths $(\sigma, h)$  there exists a unique v.h.h.
$H$ with the base $(\sigma, h)$. We say that
$\widetilde \sigma: = H|_{I_1\times \{1\}}$ is the result of the
{\it transfer of the path $\sigma$ along $h$
	with respect to the Ehresmann connection} $\mathfrak M$. 

By analogy we define v.h.h. and transfers in cases, when $I_1$ and $I_2$ are replaced
with half-intervals.

\subsection{The structure of leaves of the induced foliation}
Now we consider any smooth foliation $(M,F)$ and its graph $G(F)$.
\begin{proposition}\label{p4} Let $\mathbb L_\alpha = p^{-1}_1(L_{\alpha})$,
	$L_\alpha=L_\alpha(x)\in F$, be any leaf of the induced foliation 
	$(G(F),\mathbb F)$. Let ${\mathcal L_\alpha} = p_2^{-1}(x)$. Then:
	
	(i) The restriction $p_1|_{\mathcal L_\alpha}$ is a regular covering map of $\mathcal L_\alpha$
	onto $L_\alpha$, and its deck transformation group is isomorphic
	to the germ holonomy group $\Gamma(L_\alpha,x)$ of the leaf $L_\alpha$ at $x$.
	
	(ii) The diagonal action of $\Psi_{\alpha} = \Gamma(L_\alpha,x)$ is defined on the product
	${\mathcal L_\alpha}\times{\mathcal L_\alpha}$ of manifolds, and the quotient manifold
	$({\mathcal L_\alpha}\times{\mathcal L_\alpha})/\Psi_{\alpha}$ is
	diffeomorphic to the leaf $\mathbb L_\alpha$ and there exists a diffeomorphism
	$\Xi_{\alpha}:\mathbb{L}_{\alpha}\to (\mathcal{L}_{\alpha}\times\mathcal{L}_{\alpha})/\Psi_{\alpha}$
	such that $\Xi_{\alpha}$ is an isomorphism of the foliations $(\mathbb{L}_{\alpha}, F^{(i)})$ and
	$((\mathcal{L}_{\alpha}\times\mathcal{L}_{\alpha})/\Psi_{\alpha}, F_{i}),\, i=1,2,$ where
	$(F_{1}, F_{2})$ is the pair of foliations of complementary dimensions covered by the product $\mathcal{L}_{\alpha}\times\mathcal{L}_{\alpha}.$
\end{proposition}
	The definitions of the graph $G(F)$ of its smooth structure and of the canonical
	projections imply that the restriction $p_{1}|_{p_{2}^{-1}(x)}:p_{2}^{-1}(x)\to L_{\alpha}(x)$
	is a regular covering map with the deck transformation group isomorphic to
	$\Gamma(L_{\alpha},x).$ Hence $(i)$ is true.
	
	For the proof of $(ii)$ we consider two simple foliations
	${F}^{(i)}:=\left\lbrace p_{i}^{-1}(x)\,|\,x\in M\right\rbrace$,  $i=1,2,$
	on $G({F})$. They induce two simple foliations on the leaf $\mathbb{L}_{\alpha}$
	which will be denoted also  by  $F^{(1)}$ and $F^{(2)}$. Let us show that $TF^{(2)}$ 
	is an Ehresmann connection for the foliation $(\mathbb{L}_{\alpha}, F ^{(1)}).$
	
	Take an admissible pair $(g,f)$ of piecewise smooth paths in the leaf $\mathbb{L}_{\alpha}$
	with the common origin at $g(0)=f(0)=z_0=(x_0,<h_0>,y_0),$
	where $g: I_1\to p_2^{-1}(y_0)$ and
	$f: I_2\to p_1^{-1}(x_0)$, $I_1 = I_2 = [0,1].$ Let $g_0=p_1\circ g.$
	Denote by $g_{t}$, $t\in I_2$, the path with the origin in $f(t)$
	covering the path $g$ relatively the covering  map $p_{1}|_{p_{2}^{-1}(y_t)}$,
	where $y_t =  p_2(f(t))$.
	Put $H(s,t):= g_t(s)\,\,\, \forall (s,t)\in I_1\times I_2.$
	It is easy to check that $H$ is a vertical-horizontal homotopy with the base
	$(g, f).$ Hence $TF^{(2)}$ is an integrable Ehresmann connection
	for $(\mathbb L_{\alpha} ,F^{(1)}).$
	
	Let $k: \widetilde{\mathbb L}\to \mathbb L_{\alpha}$ be a universal covering map for the leaf $\mathbb{L}_{\alpha}$
	and $\widetilde{F^{(i)}}= k^{*}F^{(i)}$, $i = 1,2,$ be an induced foliations on $\widetilde{\mathbb L}$.
	It is easy to show that $T\widetilde{F^{(2)}}$ is an integ\-rable Ehresmann connection for
	the foliation $\widetilde{F^{(1)}}$. Due to simple connectivity of $\widetilde{\mathbb{L}}$
	according to the decomposition theorem proved by S. Kashiwabara \cite{Kash}, we
	identify $\widetilde{\mathbb L}$ with the product manifold $\widetilde{N_1}\times\widetilde{N_2}$.
	In this case $\widetilde{F^{(1)}} = \{\widetilde{N_1}\times\{v\}\,|\, v\in\widetilde{N_2}\}$,
	$\widetilde{F^{(2)}} = \{\{u\}\times\widetilde{N_1}\,|\, u\in\widetilde{N_1}\}$.
	Thus, $(\mathbb L_\alpha,F^{(1)},F^{(2)})$ is a pair of two transverse
	simple foliations on $\mathbb L_\alpha$ covered by the product of manifolds.
	
	In the work \cite{ShZh} that construction is said to be a {\it
		simple transversal bifibration} and is denoted by $(\mathbb
	L_\alpha,p_1,p_2,L_\alpha,L_\alpha)$.
	
	Let us fix a point $z = (x, <1_{x}>,x)\in\mathbb{L}_{\alpha}.$
	According $(i)$ the group $\Psi_{\alpha}$ acts on the covering
	manifold $\mathcal{L}_{\alpha}$ as a group of deck transformations and the diagonal
	action of $\Psi_{\alpha}$ on the product manifolds $\mathcal{L}_{\alpha}\times\mathcal{L}_{\alpha}$
	is defined. This action is proper and totally discontinuous, and preserves the product.
	The quotient manifold $(\mathcal{L}_{\alpha}\times\mathcal{L}_{\alpha})/\Psi_{\alpha}$
	with two foliations $F_{1}, F_{2}$ whose leaves are covered by the respective leaves of the 
	trivial foliations of the product $\mathcal{L}_{\alpha}\times\mathcal{L}_{\alpha}$ is defined.
	
	Proposition 1 from \cite{ShZh} implies the existence of a diffeomorphism
	$$\Xi_{\alpha}:\mathbb{L}\to (\mathcal{L}_{\alpha}\times\mathcal{L}_{\alpha})/\Psi_{\alpha}$$ such that
	$\Xi_{\alpha}$ is an isomorphism in the foliation category $\mathfrak F\mathfrak o\mathfrak l$ of
	both pairs of foliations $(\mathbb{L}_{\alpha}, F^{(i)})$ and
	$((\mathcal{L}_{\alpha}\times\mathcal{L}_{\alpha})/\Psi_{\alpha}, F_{i}),\, i=1,2.$

\subsection{Proof of Theorem \ref{t2}}

Denote by $(M, F)$ an arbitrary pseudo-Riemannian foliation on a
pseudo-Riemannian manifold $(M,g)$.

Observe that a Riemannian structure may be considered as a $G$-structure of the 
first order \cite{Ko}. As for as every pseudogroup of local isometries of a 
pseudo-Rieman\-ni\-an  manifold is quasi-analytic we may apply Winkelnkemper's 
criterion indicated above, according to which the graph $G({F})$ is Hausdorff.

Observe that the  definition of the inducted foliation $(G(F), \mathbb F)$ implies
that both foliations $(M, F)$ and $(G(F), \mathbb F)$ are given by the same
holonomy pseudogroup $\mathcal H$. Recall that the germ holonomy group of an arbitrary
foliation is interpreted as a group of germs at the relevant point $v$ of the local
transformations $\varphi$ from the holonomy pseudogroup $\mathcal H$ such that $\varphi(v) = v$
(see, for example, \cite{Sal}).

This implies that the holonomy groups $\Gamma (L,x)$ and $\Gamma (\mathbb{L},z)$,
where $\mathbb{L} = p_1^{-1}(L)$, $z =(x,<h>,y)\in G(F)$, are isomorphic.

Thus the assertion 1 of Theorem \ref{t2} is valid.

Let $\mathbb F$ be the induced foliation and
${F}^{(i)}:=\left\lbrace p_{i}^{-1}(x)\,|\,x\in M\right\rbrace, i=1,2,$
be two simple foliations on the graph $G(F)$.
Define a special pseudo-Riemannian metric on the graph $G(F)$.
Denote by $\mathfrak{M}$ the $q$-dimensional distribution on $M$ orthogonal to the pseudo-Riemannian
foliation $(M,F).$ Since the pseudo-Riemannian metric $g$ is non-degenerate on the leaves
of this foliation, there exists a decomposition of the tangent space $T_{x}M,\,\,x\in M,$ of $M$
into the orthogonal sum of vector subspaces: $$T_xM = T_x F\oplus\mathfrak M_x.$$
For any $z =(x,<h>,y)\in G(F)$ put
$$\mathfrak{N}_z :=\{ X\in T_{z}G(F)\,|\,p_{1*z}X\in\mathfrak{M}_x, \, p_{2*z}(X)\in\mathfrak{M}_{y}\}.$$
Emphasize that there exists a bijective mapping of the a intersection $p_{1}^{-1}(x)\cap p_{2}^{-1}(y)$ to the holonomy group $\Gamma(L, x)$. Therefore
$\mathfrak{N} = \{\mathfrak{N}_{z}\,|\,z\in G(F)\}$ is a smooth $q$-dimensional
distribution on the graph $G(F)$ and, for any $z\in G(F),$  the tangent vector space
$T_{z} G(F),$ admits the following decomposition into a direct sum of vector subspaces
\begin{equation}\label{q0}
T_{z}(G(F)) = T_{z}(F^{(1)})\oplus\mathfrak{N}_{z}\oplus T_{z}( F^{(2)}).
\end{equation} According to the decomposition (\ref{q0}),
any vector field $X\in\mathfrak{X}(G(F))$ may be represented in the form
\begin{equation}\label{q1}
X = X^{(1)} + X^{\mathfrak{N}} + X^{(2)},
\end{equation}
where $X^{(i)}\in\mathfrak{X}_{F^{(i)}}G$ and $X^{\mathfrak{N}}\in\mathfrak{X}_{\mathfrak{N}}G(F).$
Let us define the pseudo-Riemannian metric $d$ on $G(F)$ by the equality
\begin{equation}\label{q2}
d(X,Y): = g(p_{1*}X, p_{1*}Y) + g(p_{2*}X^{(1)}, p_{2*}Y^{(1)}),
\end{equation}
where $X, Y$ are represented in the form $(\ref{q1}).$

Observe that $$d(X,Y): = g(p_{1*}X^{(2)} ,p_{1*}Y^{(2)})+g(p_{1*}X^{\mathfrak{N}},
p_{1*}Y^{\mathfrak{N}}) + g(p_{2*}X^{(1)} ,p_{2*}Y^{(1)}), $$
for any  $X, Y\in \mathfrak{X}(G).$ Note that the foliations $F^{(1)},F^{(2)}$ 
and distribution $\mathfrak{N}$ are pairwise orthogonal in the pseudo-Riemannian manifold $(G(F),d).$ 
Moreover, the restriction of $d$ onto any leaf of each foliation $\mathbb{F}, F^{(1)}$ 
and $ F^{(2)}$ on the graph $G(F)$ is non-degenerated, and, for any $X, Y\in \mathfrak{X}_{T\mathbb F}(G):$

$$d(X,Y): = g(p_{1*}X^{(2)}, p_{1*}Y^{(2)})+ g(p_{2*}X^{(1)}, p_{2*}Y^{(1)}).$$
This means that the induced pseudo-Riemannian metric $d |_{\mathbb{L}}$ on a leaf $\mathbb{L}$ 
is locally a direct product of the 
pseudo-Riemannian metric induced on leaves of foliations $F^{(1)}|_{\mathbb L}$ and
$F^{(2)}|_{\mathbb L}$. Therefore (see, for example \cite{Wu} or \cite{KN} for the Riemannian case)
the distributions  $TF^{(1)}$ and $TF^{(2)}$ are orthogonal and parallel on $(\mathbb{L},d |_{\mathbb{L}})$.
Hence $(\mathbb{L}, d|_{\mathbb{L}})$ is a non-degenerate reducible pseudo-Riemannian manifold.

Further the restriction of $d$ (respectively, $g$ or $p_{i*}$) onto the vector subspaces
of corresponding vector subspaces of $T_{z}G(F)$, $z\in G(F)$, will be also denoted also $d$
(respectively, $g$ or $p_{i*}$).

The canonical projection $p_{1}:G(F)\to M$ is pseudo-Riemannian because
$p_{1*z}:T_{z}F^{(2)}\oplus\mathfrak{N}\to T_{x}M$ is an isomorphism of
the pseudo-Euclidean vector spaces $(T_{z}F^{(2)}\oplus\mathfrak{N}_{z},d)$ and
$(T_{x}F,g|_{T_{x}F})$ by definition of $d.$

Let us show that the canonical projection $p_{2}$ is also a pseudo-Riemannian submersion.
Take any point $z=(x,<h>,y)$ in $G(F).$ The pseudo-Rieman\-ni\-ance of
$(M,F)$ implies the existence of a linear isomorphism
$\phi_{xy}:\mathfrak{M}_{x}\to \mathfrak{M}_{y}$ induced by the local holonomy isometry along $h$ of the pseudo-Euclidean
vector spaces $(\mathfrak{M}_{x},g)$ and $(\mathfrak{M}_{y},g).$
According to the definition of $d$ the restriction $p_{2*z}:T_{z}F^{(1)}\to T_{y}F$
is an isomorphism of pseudo-Euclidean spaces.
Observe that the restriction
$p_{2*z}:T_{z}F^{(1)}\oplus\mathfrak{N}_{z}\to T_{y}M = T_{y}F\oplus\mathfrak{M}_{y}$
is equal to $p_{2*z} = (\phi_{xy}\circ p_{1*z}|_{\mathfrak{N}_{z}}, p_{2*z}|_{T_{z}F^{(1)}}).$
This implies that $p_{2}$ is also a pseudo-Riemannian submersion.

Take any vector $X\in \mathfrak{N}_{z}$. Let $\gamma=\gamma(s)$ be 
the geodesic passing through the point $z$ in the direction of the vector $X$, 
i.e. $\gamma(0)=z,\,\dot{\gamma}(0)=X.$ Therefore $\gamma$ is the geodesic orthogonal 
to the leaves $L^{(1)}(z)$ and $L^{(2)}(z)$ of the foliations $F^{(1)}$ and $F^{(2)}$. 
As for as $L^{(1)}$ and $L^{(2)}$ are the fibers of the pseudo-Riemannian submersions 
$p_{1}$ and $p_{2}$, according to Corollary $1$ at the any point of $\gamma(s)$ the tangent 
vector $\dot{\gamma}(s)$ is orthogonal to the both fibers $L^{(1)}(\gamma(s))$ and $L^{(2)}(\gamma(s))$.  
Since $T_{\gamma(s)}\mathbb{L}=T_{\gamma(s)}L^{(1)}\oplus T_{\gamma(s)}L^{(2)},$ then the tangent vector 
$\dot{\gamma}(s)$ is orthogonal to $T_{\gamma(s)}\mathbb{L}.$ Thus geodesic $\gamma=\gamma(s)$ 
of the pseudo-Riemannian manifold $(G(F),d)$ orthogonal to leaves of the foliation $(G(F),\mathbb{F})$ 
at the one its point is orthogonal to the foliation $(G(F),\mathbb{F})$ at each its point.

According to Theorem $1$, due to non-degeneracy of a pseudo-Riemannian metric on leaves 
of this foliation, $(G(F),\mathbb{F})$ is a pseudo-Riemannian foliation. 

Assume that there exists another pseudo-Riemannian metric $\hat d$ on $G(F)$
satisfying the second statement of Theorem \ref{t2}. Let $\hat{\mathfrak N}$ be the $q$-dimensional
distribution that is orthogonal to the foliation $\mathbb F$ in $(G(F),\hat d)$.
Therefore $p_{1*}\hat{\mathfrak N}_z = {\mathfrak M}_x$ and $p_{2*}\hat{\mathfrak N}_z = {\mathfrak M}_y$
for every point $z=(x,<h>,y)$ in $G(F).$
Consequently $\hat{\mathfrak N} = \mathfrak N$ and $\hat d(X,Y) = d(X,Y) = g(p_{1*}X,p_{1*}Y)$
for any $X,Y\in\mathfrak X_{\mathfrak N}G(F)$. According to our assumption,
in relation to both metrics $d$ and $\hat d$ the foliations $F^{(1)}$, $F^{(2)}$
are orthogonal, with $p_i$, $i=1,2$, are pseudo-Riemannian submersions.
Due to the decomposition (\ref{q1})
and bilinearity of $d$ and $\hat d$ it is necessary that $d = \hat d$.

Thus the statements 1 and 2 of Theorem~\ref{t2} are proven.
The statement 3 of Theorem~\ref{t2} follows from Proposition~\ref{p4}.

\section{Two graphs of a foliation with an Ehresmann connection}\label{sec:4}
\subsection{Holonomy groups of foliations with Ehresmann connections.}

Let $(M, F)$ be a foliation with an Ehresmann connection $\mathfrak M$.
Take any point $x\in M$. Denote by $\Omega_x$
the set of horizontal curves with the origin at $x$.
An action of the fundamental group $\pi_1(L,x)$ of the leaf $L=L(x)$
on the set $\Omega_x$ is defined in the following way:
$$\Phi_x:\pi_1(L,x)\times \Omega_x\to \Omega_x:([h],\sigma)\mapsto\widetilde\sigma,$$
where $[h]\in\pi_1(L,x)$ and $\widetilde\sigma$ is the result of
transfer of $\sigma\in\Omega_x$ along $h$ in relation to $\mathfrak M$.
Let $K_{\mathfrak M}(L,x)$ be the kernel of the action $\Phi_x$, i.e.
$$K_{\mathfrak M}(L,x)=\{\alpha\in \pi_1(L,x)\,|\, \alpha(\sigma)=\sigma \,\, \, \forall \sigma\in \Omega_x\}.$$
The quotient group $H_{\mathfrak M}(L,x)=\pi_1(L,x)/K_{\mathfrak M}(L,x)$
is the {\it $\mathfrak M$-holonomy group of the leaf} $L$, see \cite{BH}. Due to the pathwise
connectedness of the leaves, the $\mathfrak M$-holonomy groups at
different points on the same leaf are isomorphic.

Let $\Gamma (L,x)$ be a germinal holonomy group of the leaf $L$. Then
there exists a unique group epimorphism $\chi :H_{\mathfrak
	M}(L,x)\to \Gamma (L,x)$ satisfying the equality
$$\chi \circ \mu=\nu , \eqno(1) $$ where
$\mu:\pi_1(L,x)\to H_{\mathfrak M}(L,x):[h]\mapsto[h]\cdot K_{\mathfrak M}(L,x)$ 
is the quotient map and
$\nu:\pi_1(L,x)\to\Gamma(L,x):[h]\mapsto<\!\!h\!\!>$, where $<\!\!h\!\!>$ 
is a germ of the holonomy diffeomorphism of
a transverse $q$-dimensional disk along the loop $h$ at the point $x$.

Let us emphasize that the $\mathfrak M$-holonomy group $H_{\mathfrak M}(L,x)$
has {\it a global character} unlike the germinal holonomy group
$\Gamma(L,x)$ having a local-global character: global along the leaves and
local along the transverse directions.

\subsection{The graph $G_{\mathfrak M}(F)$}
The graph $G_{\mathfrak M}(F)$ of a foliation $(M, F)$ with an Ehresmann connection
was introduced by the first author in \cite{ZhG} (see also \cite{ZhL}).

Let $(M,F)$ be a foliation of arbitrary dimension $k$
on a $n$-manifold $M$ and $q=n-k$ be the codimension of this
foliation. Suppose that the foliation $(M, F)$ admits an
Ehresmann connection $\mathfrak M.$

Consider the set $\Omega_x$ of $\mathfrak M$-horizontal curves with the origin at $x\in M$.

Take any points $x$ and $y$ in the leaf $L$ of $(M, F)$.
Introduce an equivalence relation $\rho$ on the set $A(x,y)$ of
vertical paths in $L$ connecting $x$ with $y$. Paths $h$ and $f$
in $A(x,y)$ are called $\rho$-equivalent if the loop $h\cdot
f^{-1}$ generates the trivial element of the $\mathfrak M$-holonomy
group $H_{\mathfrak M}(L,x)$. In other words, paths $h$ and $f$ are
$\rho$-equivalent iff they define the same transfers of $\mathfrak
M$-horizontal curves from $\Omega_x$ to $\Omega_y$ in relation to the Ehresmann
connection $\mathfrak M$. The $\rho$-equivalence class containing $h$ is
denoted by $\{h\}$.

The set of ordered triplets $(x,\{h\},y)$, where $x$ and $y$ are
any points of an arbitrary leaf $L$ of the foliation $(M,F)$ and $\{h\}$
is a class of $\rho$-equivalent paths from $x$
to $y$ in $L$, is called {\it the graph of the foliation
	$(M,F)$ with an Ehresmann connection $\mathfrak M$} and is
denoted by $G_{\mathfrak M}(F)$.

The maps
$$p_1:G_{\mathfrak M}(F)\to M: (x,\{h\},y)\mapsto x,\,\,
p_2:G_{\mathfrak M}(F)\to M: (x,\{h\},y)\mapsto y$$ are called the
{\it canonical projections}. 

The graph $G_{\mathfrak M}(F)$ is
equipped with a smooth structure and the binary operation
$(y,\{h_1\},z)\ast(x,\{h_2\},y):=(x,\{h_1\cdot h_2\},z)$ becomes a smooth $\mathfrak M$-holonomy groupoid.

We summarize the first author's results (\cite{ZhG} and  \cite{ZhL}) into the following theorem.

\begin{theorem}\label{t4} 1. The graph $G_{\mathfrak{M}}(F)$ of
	a foliation $(M, F)$ with an Ehresmann connection
	$\mathfrak M$ admits a structure of a Hausdorff manifold and the canonical
	projections $p_1$ and $p_2$ determine locally trivial fibrations with common typical fibre $Y$,
	and for any point $x\in M$ there is a regular covering map
	$p_x: Y\to L(x)$ with the deck transformation group isomorphic to
	the $\mathfrak{M}$-holonomy group $H_\mathfrak{M}(L,x)$.
	
	2. The distribution
	$\mathfrak N: = \{\mathfrak N_z \,|\, z=(x,\{h\},y)\in G_\mathfrak{M}(F)\}$ where
	$$\mathfrak N_z = \{X\in T_z G_\mathfrak{M}(F) \,|\, p_{1*z}X\in\mathfrak M_x,  p_{2*z}X\in\mathfrak M_y\}$$
	is an Ehresmann connection for the induced foliation
	$\mathbb F:=\{p_1^{-1}(L_{\alpha})\,|\,L_{\alpha}\in F\}.$
	
	3. The group epimorphism
	$\chi: H_{\mathfrak M}(L,x)\to \Gamma (L,x)$ satisfying the equality
	$(1)$ is a group isomorphism
	and the map $$G_{\mathfrak M}(F)\to G(F):
	(x,\{h\},y)\to (x,<h>,y)$$
	is a groupoid isomorphism  if and only if the holonomy
	pseudogroup ${\mathcal H}(M, F)$ of the foliation $(M,F)$ is quasi-analytic.
\end{theorem}

\section{The proof of Theorem \ref{t3}}\label{sec:5}

I. Let us prove the following lemma which will be useful in this proof.

\begin{lemma}\label{l1}
	Let $(M, F)$ be any transversally complete pseudo-Riemannian foliation.
	Then the orthogonal $q$-dimension distribution $\mathfrak{M}$ is an Ehresmann
	connection for $(M, F).$
\end{lemma}
	Theorem~\ref{t1} implies that for any point $x\in M$ there exists an adapted
	chart $(U,\varphi)$ such that projection $p_{U}: U\to V=U/F_{U}$ is a pseudo-Riemannian
	submersion, and $\mathfrak{M}_{U}$ is a geodesically invariant distribution. Hence there is
	a chart $(\tilde{U},\tilde{\varphi})$ at $x\in\tilde{U}\subset U,\,\tilde{\varphi}=\varphi|_{\tilde{U}}$
	with the following property: for any admissible pair of the paths $(\gamma, f)$,
	with origin $\gamma(0)=f(0) = x$ in $\tilde{U}$, there exists a vertical-horizontal homotopy
	$H$ with the base $(\gamma, f)$ where $\gamma$ is a horizontal geodesic. 
	The application of the results from \cite{ZhP} to transversally complete pseudo-Riemannian
	foliations implies that all statements for foliations with an Ehresmann connection remain 
	true if we prove the existence of a vertical-horizontal homotopy only for a pair
	of the form $(\gamma, f),$ where $\gamma$ is a horizontal geodesic.
	
	Let us take any admissible pair of paths $(\gamma, f)$ where $\gamma$ is a
	horizontal geodesic. Cover the path $f(I_{2})$ by a finite chain of adapted charts
	$\tilde{U}_{1},\tilde{U}_{2},...,\tilde{U}_{k}$ which have the property indicated above.
	Therefore there is a number $\delta >0$ such that for the pair $(\gamma|_{[0,\delta]}, f)$
	there exists the v.h.h. $H_{\delta}:[0,\delta]\times I_{2}\to\tilde{U}_{1}\cup\tilde{U}_{2}\cup...\cup\tilde{U}_{k}$
	with the base $(\gamma|_{[0,\delta]},f)$.
	We continue this construction and get the prolongation $H$ onto the
	$[0,\delta_{0})\times I_{2},$ where $\delta_{0}\in (\delta,1),$ and $H$ is the
	vertical-horizontal homotopy with the base $(\gamma|_{[0,\delta_{0})},f).$
	Hence for every $t\in I_{2}$ we have the transfer
	$\gamma_{f(t)}(s): = H(s,t),\, s\in [0,\delta_{0}).$ Due to transversal
	completeness of the foliation $(M, F)$, the geodesic $\gamma_{f(t)}(s)$ is
	defined for any $s\in I_{1}.$ Therefore we may put
	$H(s,t):=\gamma_{f(t)}(s)\,\,\forall (s,t)\in I_{1}\times I_{2}$.
	The verification shows that $H$ is the vertical-horizontal homotopy with the base
	$(\gamma, f)$,  that finishes the proof of Lemma~\ref{l1}.

Let $(M, F)$ be a transversally complete pseudo-Riemannian foliation.
According to Lemma \ref{l1}, the distribution $\mathfrak{M}$ is an Ehresmann connection for
$(M,F).$ Therefore two graphs $G(F)$ and $G_{\mathfrak{M}}(F)$ are defined.
Since the holonomy pseudogroup ${\mathcal H}(M, F)$ of this foliation is quasi-analytic,
according to the statement~3 of Theorem~\ref{t4} the holonomy group $\Gamma(L, x)$ 
and $H_{\mathfrak{M}}(L, x)$ are isomorphic for any $x\in M$ and we may identify 
the graph $G(F)$ with the graph $G_{\mathfrak M}(F)$.

In this case in accordance with the first statement of Theorem~\ref{t4}
the canonical projections  $p_{i}: G(F)\to M$, $i= 1,2$, define locally trivial bundles
with the same typical fibre $Y$. It is easy to see that $Y = L_0$ is diffeomorphic to any leaf
$L_\alpha$ without the holonomy of the foliation $(M, F)$.

By statement~2 of Theorem~\ref{t4} the $q$-dimensional
distribution $\mathfrak{N}$ is an Ehresmann connection for the induced
foliation $(G(F),\mathbb{F})$ on the manifold $(G(F),d)$.
It implies that the distribution $\mathfrak{K}=\mathfrak N \oplus TF^{(2)}$ is an Ehresmann connection for the
submersion $p_1: G(F)\to M.$ Moreover, a $\mathfrak{K}$-lift of any $\mathfrak M$-curve in $M$
is a $\mathfrak N$-curve in $G(F)$. Since $p_1$ is a pseudo-Riemannian
submersion, any $\mathfrak N$-lift of a $\mathfrak M$-geodesic $\gamma$ in $(M,g)$ to
any point $z\in p_1^{-1}(\gamma(0))$ is a $\mathfrak N$-geodesic in $(G(F),d)$
and the projection $p_1\circ\hat\gamma$ of every $\mathfrak N$-geodesic
$\hat\gamma$ in $(G(F),d)$ is the $\mathfrak M$-geodesic in $(M,g)$. These facts
and the transversal completeness of $(M,F)$ imply the transversal
completeness of $(G(F),\mathbb F)$.

Since the holonomy pseudogroup ${\mathcal H}(G(F), \mathbb F)$ is also quasi-analytical,
according to statement~3 of Theorem~\ref{t4} the holonomy groups
$\Gamma(\mathbb{L}_\alpha, z)$ and $H_{\mathfrak{N}}(\mathbb{L}_\alpha, z)$ are isomorphic.
In accordance with the first statement of Theorem \ref{t2}
the holonomy groups  $\Gamma(\mathbb{L}{_\alpha},z)$ and $\Gamma(L_\alpha,x)$,
$x=p_1(z)$, are isomorphic.
Therefore $\Gamma(L_\alpha)\cong H_{\mathfrak M}(L_\alpha)\cong
\Gamma(\mathbb{L}_{\alpha})\cong H_{\mathfrak N}(\mathbb{L}_{\alpha})$.

Thus, three statements of Theorem \ref{t3} are proved.

Remark that statement 4 follows from the proven statements 2 and
3 of Theorem~\ref{t2}.

II. The following characteristic property is well known for the
geodesically invariant foliations $(M, F)$ on Riemannian manifolds $M$ and by analogy
is proved for pseudo-Riemannian manifolds $M$ (see, for example \cite {Y}, Proposition 2.7).
Denote by $\mathcal{L}_{X}$ the Lie derivative along the vector field $X.$

\begin{proposition}\label{pL}
	Let $(M, g)$ be pseudo-Riemannian manifold and $F$ be a codimension $q$ foliation of $M.$
	Then $F$ is  geodesically invariant if and only if $(\mathcal{L}_{X}g)(Y,Z)=0$ for all
	$X\in\mathfrak{X}_{\mathfrak{M}}(M)$ and for all $Y, Z\in\mathfrak{X}_{F}(M),$ where
	$\mathfrak{M}$ is the distribution which consists of all vectors orthogonal to $TF.$
	
	Furthermore if $F$ is geodesically invariant and $X\in\mathfrak{X}_{\mathfrak{M}}(M)$
	is a local foliated vector field which can define a local $1$-parameter group,
	then a local $1$-parameter group of local transformations generated by
	$X$ preserves the metrics induced on plaques.
\end{proposition}

It is known that for any two leaves $L_{0}$ and $L$ there exists a
piecewise smooth horizontal geodesic $\sigma:[0, 1]\to M$ such that $a_{0}=\sigma(0)\in L_{0}$
and $a=\sigma(1)\in L.$ Let $L_{0}$ be a fixed leaf without holonomy and $L$
be any other leaf of the foliation $(M, F).$

Let $x$ be an arbitrary point in $L_{0}$. Connect $a_{0}$ with $x$ by a vertical path
$h:[0,1]\to L_{0},\,\,\,h(0)=a_{0}, \,\,h(1)=x.$ As $(M, F)$ is a transversally
complete pseudo-Riemannian foliation according to Lemma \ref{l1}
$\mathfrak{M}$ is an Ehresmann connection for foliation $(M, F).$ Then for any
admissible pair $(\sigma, h)$ there exists a vertical-horizontal homotopy
$H: I_{1}\times I_{2}\to M$ with the base $(\sigma, h).$
Let $\tilde{\sigma}:= H|_{I_{1}\times\{1\}}.$ In this case [\cite{ZhL}, Lemma 1]
the map $f_{\sigma}: L_{0}\to L: x\mapsto\tilde{\sigma}(1)$
is well defined and is a regular covering with the deck transformation group isomorphic
to $H_{\mathfrak M}(L, x) \cong \Gamma(L, x).$

According to our assumption, the foliation $(M, F)$ is geodesically invariant, then
Propositi\-on~\ref{pL} implies that the covering $f_{\sigma}:L_{0}\to L$ is local isometry
in relation to the induced pseudo-Riemannian metric $g|_{L_{0}}$ and $g|_{L}$ on $L_{0}$ and $L.$

Therefore for any leaf $L$ without holonomy  the map $f_{\sigma}: L_{0}\to L$ is an isometry.

Using Proposition \ref{pL} and considering that $p_{1}$ and $p_{2}$
are pseudo-Rieman\-ni\-an submersions we get that $\mathbb{F}$, $F^{(1)}$ and
$F^{(2)}$ are geodesically invariant and pseudo-Riemannian foliations on
the pseudo-Riemannian manifold $(G(F),d)$, i.e. $(3)$ is true.

Now it easy to check the fulfilment of statements $(ii)$ and $(iii)$ in II.

\section{Lorentzian foliations of codimension 2 on closed 3-manifolds}\label{sec:6}
\subsection{Theorem of C. Boubel, P. Mounoud, C. Tarquini \cite{BMT}}
\begin{definition} The algebraic Anosov flows, up to finite coverings and finite quotients, are  the following:
	\begin{enumerate}
		\item The geodesic flows of the unit tangent bundle of hyperbolic compact surfaces.
		\item The flows defined by the suspensions of linear hyperbolic diffeomorphisms of the $2$-torus.
	\end{enumerate}
\end{definition}

An application of Molino's theory of Riemannian foliations on compact ma\-nifolds \cite{Mo}
and the classification of the Lorentzian Anosov flows given by  E.~Ghys in \cite{Gh}
allowed C. Boubel, P. Mounoud, C. Tarquini [\cite{BMT}, {Theorem 4.1}] to describe
the structure of transversally complete  Lorentzian foliations of codimension~2
on closed 3-manifolds in the following way:
\begin{theorem}
	Up to finite coverings, a 1-dimensional transversally complete Lorentzian
	foliation on a compact closed 3-manifold is either smoothly equivalent to a foliation
	generated by an algebraic Anosov flow or a Riemannian foliation.
\end{theorem}

The structure of suspended algebraic Lorentzian foliations of codimension 2
on closed 3-manifolds and their graphs is described in the following subsection \ref{E1}.

\subsection{Example}\label{E1}
Let $A=\begin{pmatrix} a & b \\ c & d \end{pmatrix}\in SL(2,\mathbb{Z})$
be a matrix of integers such that $ad-bc=1$ and $a+d > 2$. Such matrix $A$
induces an Anosov automorphism $f_{A}$ of a torus $\mathbb{T}^2:=\mathbb{R}^2/\mathbb{Z}^{2}$ conserving its orientation.

Let us consider the action of the group of integers $\mathbb{Z}$ by the formula
\begin{equation}
\Phi_{A}:=\mathbb{T}^2\times\mathbb{R}^1\times\mathbb{Z}\to\mathbb{T}^2\times\mathbb{R}^1:(u,t,n)
\mapsto (f_{A}^{n}(u),f+n), \,\, n\in\mathbb{Z}.
\end{equation}
Then the quotient manifold $M:=\mathbb{T}^{2}\times_{\mathbb{Z}}\mathbb{R}^1$ is defined, and $M$ is a closed $3$-manifold. Let $\varphi:\mathbb{T}^{2}\times\mathbb{R}\to M$ be the quotient mapping. Since $\Phi_{A}$ keeps the trivial foliation $F_{tr}=\{ \{u\}\times\mathbb{R}^1|u\in\mathbb{T}^2\}$, the foliation $(M, F)$ of codimension $2$ is defined.

There exists the Lorentzian  metric $g = \eta\begin{pmatrix}
-2c & a-d \\
a-d & 2b
\end{pmatrix}$, $\eta\in R\setminus\{0\},$
on the plane $\mathbb{R}^2$ which is invariant in relation to $A$ (\cite{ZhR}, Theorem 1). This metric induces the flat Lorentzian metric on $\mathbb{T}^2$ which is denoted also by $g.$ Therefore $\tilde{g}:=\begin{pmatrix}
-2c & a-d & 0 \\
a-d & 2b & 0 \\
0 & 0 & 1
\end{pmatrix} $ is a flat Lorentzian metric on the manifold $M,$ and $(M, \tilde{g})$
is a locally pseudo-Euclidean manifold with the Lorentzian totally geodesic foliation $(M, F).$

Let $k: \tilde{M}\to M$ be the universal covering map for $M$ and $\tilde{F}:=k^* F$
be the induced foliation $\tilde{M}.$ Then
$\tilde{M}=\mathbb{R}^{3}\cong\mathbb{R}^{2}\times\mathbb{R}^1$ and
$\tilde{F} = \{\{v\}\times\mathbb{R}^{1}\,|\,v\in\mathbb{R}^{2}\}.$ Let
$pr:\tilde{M}=\mathbb{R}^{2}\times\mathbb{R}^1\to\mathbb{R}^{2}$
be the projection onto the first multiple.

The group $\Psi:= <A>\cong\mathbb{Z}$ is the global holonomy group of
the foliation $(M, F)$ that is covered by the trivial fibre bundle
$pr:\tilde{M}=\mathbb{R}^{2}\times\mathbb{R}^1\to\mathbb{R}^{2}$.
Emphasize that $(M, F)$ is Riemannian foliation if and only if there exists a Riemannian metric $k$ on $\mathbb{R}^2$ such that $\Psi$ becomes an isometry group of $(\mathbb{R}^{2}, k)$. In this case every stationary subgroup of $\Psi$ must be relatively compact. That is impossible because $\Psi$ considered as a stationary subgroup of $O_{2}\in\mathbb{R}^2$ is non-compact discrete subgroup in $Diff(\mathbb{R}^{2}).$ 

The restriction $k|_{\tilde{L}}:\tilde{L}\to L$ onto an arbitrary leaf
$\tilde{L}\cong\mathbb{R}^1$ of $(\tilde{M},\tilde{F})$ is the holonomy
covering map and a local isometry.

The graph $G(F)$ is a $4$-dimensional manifold with the induced foliation
$\mathbb{F}$ of codimensional $2$. The generic leaf $L$ of $(M, F)$
has a trivial holonomy group and $L\cong\mathbb{R}^1.$ Hence the generic
leaf $\mathbb{L}\cong L\times L \cong\mathbb{R}^2.$ The holonomy group
of any other leaf $L_{\alpha}$ of $(M, F)$ is isomorphic to $\mathbb{Z}$
and $L_{\alpha}\cong S^{1}.$ In this case the leaf $\mathbb{L}_{\alpha}: = p_{1}^{-1}(L_{\alpha})$
is isometric to the Euclidean cylinder $\mathbb{R}^1 \times_{\mathbb{Z}}\mathbb{R}^{1}.$

Let $\tilde{k}:\tilde{G}\to G(F)$ be the universal covering map, in this case
$\tilde{G}$ is diffeomorphic to $\mathbb{R}^4$ and it is provided by
the induced foliation
$\tilde{\mathbb{F}}:=\{\{u\}\times\mathbb{R}^2\,|\, u\in\mathbb{R}^2\}$,
where $\tilde{G} = \tilde{G}(\tilde{F})$ is the graph of the foliation $(\tilde{M}, \tilde{F})$
and the following diagram
\begin{center}
	$\begin{array}{ccccc}
	G(F) & \stackrel{\tilde{k}}{\longleftarrow} & \tilde{G}(\tilde{F})\cong\mathbb{R}^2\times\mathbb{R}^2 &
	\stackrel{\tilde{pr}}{\longrightarrow} & \mathbb{R}^2 \\
	\downarrow \lefteqn{p_1}& • &  \downarrow\lefteqn{\tilde{p}_1} & • &\downarrow\lefteqn{id} \\
	M & \stackrel{k}{\longleftarrow} & \tilde{M}\cong\mathbb{R}^2\times\mathbb{R}^1 & \stackrel{pr}{\longrightarrow} & \mathbb{R}^2
	\end{array}$
\end{center}\
is commutative, where $p_1: G(F)\to M$ and $\tilde{p}_1: G(\tilde{F})\to \tilde{M}$ are the canonical projections,
$\tilde{pr}:\mathbb{R}^2 \times\mathbb{R}^2 \to\mathbb{R}^2$ is the projection onto the first multiplier.

Denote by $f: \mathbb{R}^3\to\mathbb T^2\times\mathbb R^1$ the universal covering map. Let 
$y_0: = f(0)$ where $0$ is zero in $\mathbb{R}^3$ and $x_0: = \varphi(y_0)$.
Compute the fundamental group $\pi_1(M, x_0)$. Remark that the regular covering map 
$\varphi: {\mathbb T^2\times\mathbb R^1\to M}$ induces a group monomorphism 
$\widehat{\varphi}: \pi_1(\mathbb T^2\times\mathbb R^1, y_0)\to\pi_1(M, x_0)$
onto a normal subgroup $N \cong\mathbb{Z}^2$ of $\pi_1(M, x_0)$, and, for fixed point $y_0$, 
the quotient group $\pi_1(M, x_0)/N$ is isomorphic to the deck transformation group $\widehat{G}\cong\mathbb{Z}$ with
a generator $\Phi_{A}|_{\mathbb T^2\times\mathbb R^1\times\{1\}}$. Let us consider
the leaf $L_0= L_0(x_0)$ which is diffeomorphic to the circle. Observe that the inclusion 
$j: L_0\to M$ induces a group monomorphism $\widehat{j}: \pi_1(L_0,x_0)\to\pi_1(M, x_0)$ onto
$H:=Im(\widehat{j})\cong\mathbb Z$, and the deck transformation group induced by $H$ is equal $\widehat G$.
Therefore [\cite{Bump}, Proposition 1.3.1.] the fundamental group $\pi_1(M, x_0)$ is the semi-direct product
$H\ltimes N\cong\mathbb{Z} \ltimes\mathbb{Z}^2$.

Emphasize that foliations $(M, F)$ and $(G(F),\mathbb{F})$ are both totally geodesic and Lorentzian,
with the manifolds $M$ and $G(F)$ are the Eilenberg--MacLane spaces of the type $K( H\ltimes N, 1),$ i.e.
$\pi_{n}(M)=\pi_{n}(G(M)) = 0 \,\, \,\,\forall n\geq 2,\,\,\pi_{1}(M)=\pi_{1}(G(F))\cong H\ltimes N.$


\end{document}